\documentclass[12pt]{amsart}
\input xy
\xyoption{all}
\usepackage{amssymb}
\usepackage{enumerate}

\newtheorem{fact}{Fact}[section]
\newtheorem{lemma}[fact]{Lemma}
\newtheorem{theorem}[fact]{Theorem}
\newtheorem{definition}[fact]{Definition}
\newtheorem{rremark}[fact]{Remark}
\newtheorem{exa}[fact]{Example}
\newtheorem{exas}[fact]{Examples}
 \newtheorem{prob}[fact]{Problem}
\newtheorem{deta}[fact]{Details}
\newtheorem{corollary}[fact]{Corollary}

\newtheorem{observation}[fact]{Observation}
\newtheorem{proposition}[fact]{Proposition}
\newenvironment{remark}{\begin{rremark} \rm}{\end{rremark}}
\newenvironment{example}{\begin{exa} \rm}{\end{exa}}
\newenvironment{examples}{\begin{exas} \rm}{\end{exas}}

    \def\sqr#1#2{{\vcenter{\hrule height.#2pt
        \hbox{\vrule width.#2pt height#1pt \kern#1pt
            \vrule width.#2pt}\hrule height.#2pt}}}
    \def\square{\mathchoice\sqr67\sqr67\sqr{2.1}6\sqr{1.5}6}

\def\qed{~\hfill$\square$}
\def\a{\longrightarrow}
\def\fo{fo}

\newcommand{\R}{{\mathbb R}}
\newcommand{\C}{{\mathbb C}}
\newcommand{\Z}{{\mathbb Z}}

\newcommand{\N}{{\mathbb N}}
\renewcommand{\P}{{\mathbb P}}
\newcommand{\p}{{\mathbb P}}
\DeclareMathOperator{\gr}{Gr}
\DeclareMathOperator{\tp}{Tp}

\DeclareMathOperator{\tpq}{tp}
\DeclareMathOperator{\codim}{codim}

\DeclareMathOperator{\stab}{Stab}

\DeclareMathOperator{\ad}{Ad}

\DeclareMathOperator{\id}{Id}
\DeclareMathOperator{\im}{Im}

\DeclareMathOperator{\Hom}{Hom}
\renewcommand{\hom}{\Hom}

\def\E{\mathcal E}
\def\Eo{\E^0}
\def\A{\mathcal A}
\def\av{\mathcal A}
\def\K{\mathcal K}
\def\ob{\mathcal O}

\def\hol{\hbox{Hol}}
\def\phi{\varphi}
\def\u#1{\underline{#1}}

\def\kss{Kazarian spectral sequence}
\def\tpol{Thom polynomial}
\def\cd{commutative diagram}
\def\iso{\cong}
\def\htpy{\simeq}
\newcommand{\smx}[4]{\bigl(\begin{smallmatrix}{#1}&{#2}\\{#3}&{#4}\end{smallmatrix}\bigr)}

 \newcommand{\twocase}[5]{{#1}\begin{cases} {#2}& \hbox{\rm if}\ \ {#3} \cr {#4}& \hbox{\rm if}\ \ {#5}\end{cases}}

\begin{document}

\title[Calculation of Thom polynomials for group actions]{Calculation of Thom polynomials and other cohomological obstructions for group actions}
\author{L\'aszl\'o M. Feh\'er}
\address{Department of Analysis, ELTE TTK, Budapest Hungary}
\email{lfeher@math-inst.hu}
\author{Rich\'ard Rim\'anyi}
\address{Department of Mathematics, Ohio State University}
\email{rimanyi@math.ohio-state.edu}

\thanks{\noindent Supported by FKFP55/2001 as well as
OTKA D29234 (first author),
        OTKA T029759 (second author)
}

\maketitle

\section{Introduction}\label{intr}

In this paper we propose a systematic study of Thom polynomials for group actions defined by M. Kazarian in \cite{kazarian}. On one hand we show that Thom polynomials are first obstructions for the existence of a section and are connected to several problems of topology, global geometry and enumerative
algebraic geometry. On the other hand we describe a way to calculate Thom polynomials: the method of {\em restriction equations}. It turned out that though the idea is quite simple the method is very powerful. We reproduced and improved earlier result in several directions: \cite{thompol} \cite{fr-quiver}, \cite{int}, \cite{ss}. However a proper introduction to the basic theorems was missing. In this paper we try to pay this debt as well as we present the connections with obstruction theory and equivariant cohomology. We give some new results and outline possible generalizations and problems.
 \smallskip

Calculating  Thom polynomials has a long history. In retrospect the first definition of characteristic classes can be considered as the first appearance (\cite{stiefel}). But it was R. Thom who initiated their study in the case of singularities of smooth maps. The major tool for calculating them was the method of resolutions (see \cite{avgl} for an account of the method and results).

Works of  V. Vassiliev \cite{vassiliev} and M. Kazarian \cite{kazarian} clarified the connection of Thom polynomials with the underlying symmetry groups. Their works also show that the so called {\em degeneracy loci formulas} in algebraic geometry are also  Thom polynomials for group actions. So in this respect the  history of Thom polynomials dates back to the Giambelli formulas. Not surprisingly the tool of calculating them in algebraic geometry was also the method of resolutions.

Based on works of A. Sz\H ucs (\cite{szucs1}) the second author introduced a different method (the method of {\em restriction equations}) to calculate Thom polynomials for singularities of smooth maps (\cite{thompol}). After reading M. Kazarian's paper \cite{kazarian} we realized that the method of restriction equations can be easily generalized to the case of Thom polynomials for group actions. Similar type of methods were used first in \cite{ab}, see in Section \ref{kss}.

The paper is organized the following way: In Section \ref{obstruction} we show how Thom polynomials fit into the more general problem of finding cohomological obstructions for the existence of a section. In particular we explain the connections among Thom polynomials, first obstructions and equivariant Poincar\'e dual.
In Section \ref{rest-equations} we prove the basic Theorems \ref{alleqs} and \ref{unicity} on calculating Thom polynomials. We also make some comments on possible generalizations for other cohomology theories.
In Section \ref{early} we review some earlier results mainly from singularity theory and algebraic geometry.
In Section \ref{gl} we discuss Thom polynomials for representations of the complex groups $GL(n)$. Results on the adjoint representation and representations of $GL(2)$ are new.
In Section \ref{porteous} we study the classical case of the Giambelli-Thom-Porteous formula in details.
In Section \ref{singularities} we show how the computations of \cite{thompol} based on the generalized Pontryagin-Thom construction fit into the framework of the present paper.
In Section \ref{applications} we list some other applications of the method of restriction equations which were published elsewhere.
In Section \ref{sec:ptp} we give a simple formula to calculate Thom polynomials for the projectivization of a linear group action. As a corollary we show a simple method to calculate the degree of degeneracy loci, generalizing results of Porteous, Harris and Tu, Fulton and others.
In Section \ref{kss} we show how the Kazarian spectral sequence (the spectral sequence induced by the codimension filtration of the orbit stratification)---which were used to define the Thom polynomials in \cite{kazarian}---can be used in finding the orbits of a representation and their stabilizer groups.

We would like to thank A. Sz\H ucs for informing the second author---in 1996---that the normal bundle of an orbit of a $G$-action on a contractible manifold reduces to the stabilizer group of the points of the orbit. This was a crucial point in our computation.

\section{Obstruction Theory for Generalized Cohomology Theories}\label{obstruction}

We are looking for cohomological obstructions of sections of fiber bundles. On fiber bundle we mean a bundle $F\to E\to M$ associated to a principal $G$-bundle $G\to P\to M$---i.e. $E=P\times_G F$. This is not really a restriction (take $G$ to be the self homeomorphism group of $F$) but we are mainly interested in the case when $G$ is a Lie group and $F$ is a $G$-manifold.

Our starting point is the following elementary observation:

\begin{observation} Sections of the fiber bundle $E=P\times_G F$ are in one to one correspondence with $G$ equivariant maps from $P$ to $F$. (Here the right action on $P$ is transformed to a left action naturally.)
\end{observation}

In other words in the category of left $G$-spaces a section $s$ of $E$ is a lift of the collapse map $pt_P: P\to pt$:
\xymatrix{&F\ar[d]^{pt_F}\\P\ar[r]_{pt_P} \ar[ru]^s&pt}

Suppose now that $h$ is any contravariant functor from the category of (left) $G$-spaces into (graded) rings (the main example is equivariant cohomology). Then a section $s$ induces a commutative diagram:
\xymatrix{&h(F)\ar[ld]_{h(s)}\\h(P) &h(pt)\ar[l]^{h(pt_P)}\ar[u]_{h(pt_F)}}

So the existence of the section $s$ implies that
\begin{equation} \label{ker} \ker h(pt_F)\subset \ker h(pt_P). \end{equation}
Using the definitions
\begin{definition}
$$h(pt)=\{\mbox{\rm universal characteristic classes}\},$$
$$\im  h(pt_P)=\{\mbox{\rm characteristic classes of $P$}\},$$
$$\mathcal O_F:=\ker h(pt_F)=\mbox{\rm obstruction ideal of $F$}.$$
\end{definition}
\noindent we can reformulate (\ref{ker}) in a more familiar form, as follows.
\begin{theorem}
If the fiber bundle $E=P\times_G F$ admits a section then all characteristic classes of $P$ from the obstruction ideal $\mathcal O_F$ vanish.
\end{theorem}

So the goal of a general theory would be the calculation of $\mathcal O_F$ for $G$-spaces $F$ and $G$-equivariant cohomology theories.

\subsection{ Obstructions in Ordinary Cohomology}
From any cohomology theory $h$ we can obtain a $G$-equivariant cohomology theory via the Borel construction:
$$ h_G(F):=h(BF)\ \mbox{\rm where}\ BF=EG\times_G F.$$
If $F$ is a vector space and the $G$-action is given by the representation $\rho:G\to GL(F)$ then $BF$ is also denoted by $E_\rho$. The Borel construction is also a functor, so $G$-equivariant maps induces maps between the Borel constructions of the corresponding spaces. Notice that $BP=M$, $Bpt=BG$ and $B(pt_P): M\to BG$ is the classifying map of $P$. (Notice that $G$-equivariant $K$-theory is {\em not} obtained by the Borel construction.)

Suppose now that $F$ is $d$-2 connected i.e. $\pi_i(F)=0$ for $i<d-1$ and suppose that  $\pi_{d-1}(F)\iso \Z$. Then the---homotopy theoretical---first obstruction $\fo(E)$ is an element of $H^d(M;\Z)$ and for the universal $F$-bundle $\fo(BF)\in H^d(BG;\Z)$. Recalling that  $\fo(E)$ is the transgression of the generator of $\pi_{d-1}(F)$ in the spectral sequence of the fibration $E\to M$ we can see that  $\fo(E)$ can also be characterized (up to sign) by the following property:

\begin{proposition} If $\fo(BF)$ is not zero then it is a generator of $\mathcal O^d_F=\mathcal O_F\cap H^d(BG;\Z)\iso \Z$ and $\mathcal O^i_F=\mathcal O_F\cap H^i(BG;\Z)=0$ for $i<d$.
\end{proposition}

\subsection{Geometric realization of the first obstruction---Poincar\'e duality}

On geometric realization of a cohomology class $\alpha$ of a manifold $Y$ we mean finding a submanifold $X$ such that the Poincar\'e dual of $X$---denoted by $[X]_Y$ or $[X]$ for short---is equal to $\alpha$. Our basic example is the Euler class:

\begin{example} \label{euler}
Let $S\to M$ be an $S^{n-1}$-bundle (or an $\R^n\setminus \{0\}$-bundle) associated to a principal $SO(n)$-bundle $P$ over $M$.
Then the first obstruction is the Euler class of $S\to M$ and can be realized as the Poincar\'e dual of the zero of a transversal section of the vector bundle $E=P\times_{SO(n)}\R^n$.
\end{example}

In general if $F$ is a $G$-manifold and we want to realize $\alpha=\fo(P\times_G F)$ for some  principal $G$-bundle $P\to M$ then we try to find an open embedding of $F$ into a contractible $G$-manifold $V$. (In most applications $F$ is given as an open subset of a vector space.) Then $P\times_G V$ has a section $s$ and we expect that
$$F^c(s):=\{m\in M: s(m)\notin F\}$$
represents $\alpha$ for a generic section. The problem is that $V\setminus F$ is usually not smooth (basically the only case when it is smooth is Example \ref{euler}) so we cannot expect $F^c(s)$ to be smooth either. So we have to extend the notion of Poincar\'e dual somewhat:

\begin{definition}
Let $X$ be a closed subset of a topological space $Y$. $X$ {\em represents} the  cohomology class $[X]\in H^d(Y;\Z)$ if for the inclusion $i:Y\setminus X \to Y$
$$\twocase{\ker H^j(i)\iso}0{j<d}\Z{j=d}$$
and $[X]$ generates $\ker H^d(i)$.
\end{definition}

The class $[X]$ is defined up to sign, we can call the two choices orientations.

\begin{examples}\mbox{}
 \begin{enumerate}[(i)]
\item $X\to Y$ is a proper embedding of smooth oriented manifolds. Thom isomorphism for the normal bundle of $X$ shows that $X$  represents the Poincar\'e dual of $X$.
\item $X$ can be triangulated and represents a homology class in $Y$ and  $Y$ has a non degenerate intersection pairing (e.g. $Y$ is a compact oriented manifold with no torsion homology).
\item  $(Y,\Theta)$ is a stratified space and $X$ is the closure of some strata. Then $X$ represents a cohomology class if it represents a cocycle in the Vassiliev complex $(E^{0,*}_1(\Theta),d_1)$ where $E^{*,*}_*(\Theta)$ is the spectral sequence induced by the codimension filtration of the  stratification $\Theta$. (See \cite{vassiliev} and  \cite{kazarian} as well as Section~\ref{kss} for details.) If the strata are even dimensional this condition is automatically satisfied.
 \end{enumerate}
\end{examples}

Now almost by definition we have the following:
\begin{theorem}
Suppose that the Lie group $G$ acts on the vector space $V$ and $F$ is a $G$-invariant open subspace of $V$. Assume moreover that $\pi_i(F)=0$ for $i<d-1$ and  $\pi_{d-1}(F)\iso \Z$ and  that  for the principal $G$-bundle $P\to M$ the cohomology class $\fo(P\times_G F)$ is not zero and $F^c(s)$ represents a cohomology class in $H^d(M;\Z)$. Then $\fo(P\times_G F)=\pm[F^c(s)]$.
\end{theorem}
\begin{proof}
$\fo(P\times_G F)\in \ker H^d (i)$ for  $i:F(s) \to M$ by the naturality property of the first obstruction.
\end{proof}

\subsection{The theory of Thom polynomials}\label{tp}

We use the definition of M. Kazarian \cite{kazarian} for  Thom polynomials of group actions: Given a representation $\rho: G \to GL(V)$ we are looking for obstruction of having a section of a $V$-bundle associated to this representation {\em avoiding} a certain orbit $\eta$ (or more generally a G-invariant subset of $V$). Of course the zero section avoids any orbit different from the zero orbit but this is pathological: we want obstructions for a {\em generic} section. In effect we want to avoid the {\em closure} of  $\eta$. The Thom polynomial of $\eta$ is an obstruction for having a section in the complement $V\setminus\bar\eta$: the cohomology class represented by $\bar\eta(s)$ for a generic section $s$.

We would like to show that the Thom polynomial is the first obstruction. By the previous section we only have to calculate some homotopy groups of $V\setminus\bar\eta$. The open subspace $V\setminus\bar\eta$ is highly connected: by a simple transversality argument
$$\pi_i(F)=0\ \mbox{\rm for}\  i<d-1\ \mbox{\rm where}\ d=\codim(\eta).$$
Suppose now that $G$ is a complex Lie group acting on the complex vector space $V$. Then $\bar\eta$ is an algebraic subvariety of $V$.

\begin{observation}  \label{connectivity}
Let $X\subset \C^N$ be a $d$ (real) codimensional complex algebraic variety. Then
$$\twocase{\pi_i(\C^N\setminus X)\iso}0{i<d-1}\Z{i=d-1.}$$
(If $d=2$ then $\pi_1$ should be replaced by $H_1$.)
\end{observation}

A similar theorem holds for real varieties with some extra condition. However the closure of an orbit of a real representation is not necessarily a real algebraic variety, the algebraic closure can contain some other orbits with the same codimension.

\begin{definition}
Suppose that $V$ is a $G$-space and $\eta$ is a $G$-invariant subspace. Then the {\em avoiding ideal} of $\eta$ is the obstruction ideal of $V\setminus\bar\eta$:
$$\mathcal A_\eta:=\mathcal O_{V\setminus\bar\eta}.$$
Suppose moreover that
$$\twocase{\mathcal A^j_\eta\iso}0{j<d}\Z{j=d,}$$
Then a generator $\tp(\eta)$ of $\mathcal A^d_\eta$ is called the {\em Thom polynomial} of $\eta$.
\end{definition}

In other words $\tp(\eta)$ is the  {\em $G$-equivariant (generalized) Poincar\'e dual} of $\bar\eta$.

\begin{remark}
This definition of the  Thom polynomial is essentially the same as Kazarian's definition in \cite{kazarian}. We hided the technical details by assuming these properties of the avoiding ideal $\mathcal A_\eta$. The reason is that for the class of representations we are mostly interested in these assumptions can be easily verified.
\end{remark}

\begin{corollary} \label{zeroes}
 If $\bar\eta$ is a $d$ codimensional $G$-invariant subvariety of the complex representation  $\rho: G \to GL(V)$ we have that
$$\twocase{\mathcal A^j_\eta\iso}0{j<d}\Z{j=d,}$$
and
$$\tp(\eta)=fo\big(EG\times_G(V\setminus\bar\eta)\big).$$
Similarly if If $\bar\eta$ is the union of $k$ subvarieties, each of them are $d$ codimensional and $G$-invariant, then
$$\twocase{\mathcal A^j_\eta\iso}0{j<d}{\Z^k}{j=d.}$$
\end{corollary}

\section{Calculation of Thom polynomials---The method of restriction equations}
\label{rest-equations}

If $\eta$ is a homogeneous $G$-space---i.e. $\eta\iso G/G_\eta$ where $G_\eta$ is the stabilizer group of a point in $\eta$---then the calculation of the  obstruction ideal $\mathcal O_\eta$ can be reduced to algebra: First notice that $H^*_G(\eta)\iso H^*_{G_\eta}(pt)$ and that the map $H^*_G(pt_\eta):H^*_G(pt)\iso H^*(BG)\to  H^*_G(\eta)\iso H^*(BG_\eta)$ is equal to $H^*(Bi)$ for $i:G_\eta \to G$. Then choose compatible maximal tori for $G$ and  $G_\eta$ and describe the map $H^*(Bi)$ in terms of Chern roots. The calculation for homogeneous $G$-spaces can be considered as the first step for calculating of the  obstruction ideal for more general $G$-spaces (see Theorem \ref{intersection}). This is usually not an easy algebraic question (see Section \ref{ideal}), but if we are only interested in Thom polynomials then for a certain class of representations we have an algorithm. This algorithm is based on the observation that $\mathcal O_\eta\subset\mathcal O_\xi$ if $\xi\subset\eta$ and the following property  of the generalized Poincar\'e dual:

\begin{lemma} For any orbit $\eta$ which admits Thom polynomial $H^*(pt_\eta)(\tp(\eta))=e(\eta)$, where $e(\eta)$ is the equivariant Euler class of the normal bundle of $\eta$.
\end{lemma}

\begin{proof}
Since $\eta \to V\setminus \partial\eta $ is a proper oriented submanifold $\eta$ has a Poincar\'e dual $[\eta]$ in  $V\setminus \partial\eta$. Therefore
$[\eta]|_\eta=e(\eta)$. On the other hand $\tp(\eta)|_{V\setminus\partial\eta}=[\eta]$ (if follows from the analogous statement for non equivariant cohomology see e.g. in \cite[app. B.3]{fulton:young}).
\end{proof}

So we get many conditions for $\tp(\eta)$:

\begin{theorem} \label{alleqs} $\tp(\eta)\in \mathcal O_\xi$ for orbits $\xi\subset (V\setminus\bar\eta)$ and $H^*(pt_\eta)(\tp(\eta))=e(\eta)$.
\end{theorem}

We call the subset of these conditions where  $\codim\xi\leq\codim\eta$  the {\em restriction equations}. These equations are easier to handle since sometimes it is difficult to decide whether an orbit $\xi$ belongs to the closure of $\eta$. And for a class of representations they uniquely define $\tp(\eta)$:
 \begin{definition}
A complex vector $G$-space $V$ satisfies the {\em Euler condition} if there are finitely many orbits of  $V$ and $e(\xi)$ is not a zero divisor for any of the orbits $\xi\in V/G$.
\end{definition}

\begin{remark}
This condition on the Euler class appears first in \cite{ab} as a necessary condition for $G$-perfectness, see Section \ref{kss} for more details.
\end{remark}
Let us denote the restriction maps  $H^*(pt_\xi)$ by $j^*_\xi$. Then we have

\begin{theorem}\label{unicity}
If $V$ satisfies the  Euler condition then the  restriction equations
$$j_{\xi}^*\tp(\eta)= \begin{cases}
e(\eta) & \hbox{\rm if  } \xi=\eta \phantom{,\ \codim \xi\leq \codim \eta} \qquad \qquad \hbox{`principal equation'}\cr
            0 & \hbox{\rm if  } \xi\neq \eta,\ \codim \xi\leq \codim \eta \qquad\qquad \hbox{`homogeneous equations'} \end{cases}$$
have a unique solution.
\end{theorem}

The proof is an improvement of the discussion \cite[Sect.6]{thompol} but we translate it to the language of the obstruction ideals. The proof is based on a repeated use of the following lemma.

\begin{lemma}
Suppose that $\eta\subset\xi\cup \eta$ is a proper inclusion of (complex) manifolds and that $e(\eta)$ is not a zero divisor. Then
$$\mathcal O_{\eta\cup\xi}=\mathcal O_\eta \cap\mathcal O_\xi.$$
\end{lemma}

\begin{proof}
Let $D\eta$ denote a tubular neighborhood of $\eta$ in $\xi\cup \eta$. Replacing $G$ with its maximal compact subgroup doesn't change the equivariant cohomology groups, so we can assume that $D\eta$ is a $G$-equivariant subset of $\xi\cup \eta$.

Now looking at the $G$-equivariant Mayer-Vietoris sequence of $\xi\cup D\eta=\xi\cup \eta$:
$$\xymatrix{
H^{*-1}_G(D\eta\setminus\eta)\ar[r]^\delta&H^*_G(\xi\cup \eta)\ar[r]& H^*_G(\xi)\oplus H^*_G(\eta)\\
                                      &H^*_G(pt)\ar[u] \ar[ur]}$$
we can see that the Lemma is equivalent to the statement that $\delta=0$. Considering now the relative exact sequence of the pair $(\xi\cup \eta,\xi)$ and comparing with the Gysin sequence of the normal bundle of $\eta$ we get the commutative diagram:
$$\xymatrix{
                                                            & H^{*-n}_G(\eta)\ar[r]^{\cup e(\eta)}   & H^*_G(\eta)\\
H^{*-1}_G(D\eta\setminus\eta)\ar[r]^a \ar[drr]^(.7){\delta} & H^*_G(D\eta,D\eta\setminus\eta)\ar[r]^b \ar[u]^{\iso} & H^*_G(D\eta)\ar[u]^{\iso}\\
H^{*-1}_G(\xi)\ar[r]\ar[u]                                  & H^*_G(\xi\cup \eta,\xi)\ar[r] \ar[u]^(.3){\iso}|\hole& H^*_G(\xi\cup \eta) \ar[u] }$$
The map $b$ is injective since $e(\eta)$ is not a zero divisor. This implies that $a$ is the zero map. But $\delta$ factors through $a$ so must be zero, too.
 \end{proof}

Iterating the lemma we get the following theorem.

\begin{theorem} \label{intersection}
If $V$ satisfies the  Euler condition and $X=\bigcup\{\xi_1,\dotsc,\xi_n\}\subset V$ is open then
$$\mathcal O_X=\bigcap_{i=1}^n \mathcal O_{\xi_i}.$$
\end{theorem}\qed

\begin{proof}[Proof of Theorem \ref{unicity}]
Let $d=\codim \eta$. By Theorem \ref{intersection}  Theorem \ref{unicity} is equivalent with the statement that $\mathcal A^d_\eta=\mathcal O^d_{V\setminus\bar\eta}$ is equal to  $\mathcal O^d_{U}$ where $U=\bigcup \{\xi:\codim \xi\leq d,\ \xi\not= \eta\}$ . So it is enough to show that the inclusion $U\subset V\setminus\bar\eta$ induces an injection in degree $d$. From  the relative cohomology exact sequence
$$\xymatrix{H^d_G(V\setminus\bar\eta,U) \ar[r]&H^d_G(V\setminus\bar\eta)\ar[r]&H^d_G(U)\ar[r]&H^{d+1}_G(V\setminus\bar\eta,U)}$$
we can see that it is enough to show that \( H^d_G(V\setminus\bar\eta,U)=0\). By excision $H^*_G(V\setminus\bar\eta,U)\iso H^*_G(V,V\setminus \Sigma)$ where $\Sigma=\bigcup \{\xi:\codim \xi> d,\ \xi\not\subset \bar\eta\}$ and by Corollary \ref{zeroes} (and by looking at the relative cohomology exact sequence of \((V,V\setminus \Sigma) \)) we get that $H^*_G(V,V\setminus \Sigma)=0$ for $*<d+1$.
\end{proof}

\begin{remark}
The proof shows that the principal equation is needed only to find the generator of a subgroup isomorphic to $\Z$, so theoretically it is not necessary. Also it looks more complicated than the homogeneous equations. Strangely enough sometimes the homogeneous equations are more difficult to deal with (for example in the case of double or Kempf-Laksov Schur polynomials \cite{ss}). And in Section \ref{porteous} we show some examples where the principal equation implies the homogeneous equations. This way contrary to what is expected the principal equation itself is enough to calculate the Thom polynomials.
\end{remark}

\subsection{The role of the avoiding ideal.} \label{ideal}
One of the main advantages of the point of view on Thom polynomials presented in this section is the initiation of the avoiding ideal $\av_\eta$. First, its definition is a very straightforward topological idea, second its meaning fits into general topological studies: the ideal is the collection of ``first obstructions'' to a section avoiding $\eta$. Although when we concentrated on computation, we gave Theorems \ref{alleqs}, \ref{unicity} only on one element of the avoiding ideal, the Thom polynomial. However, our proof also says, that the elements of the avoiding ideal are {\em exactly} the solutions of the homogeneous equations. Thus we obtained a description of the avoiding ideal as a kernel of a homomorphism---which is in most situations the desired description. In some other cases one wants to have a generator set for this ideal, this needs additional algebraic work. This was done in the classical situation by P. Pragacz \cite{pragacz-enum}, see also Proposition \ref{pr-ideal}. \smallskip

The greatest importance of the notion of $\av_\eta$ is, however, that this is the notion which generalizes to extraordinary cohomology theories, because of the simplicity of its definition. The notion of Thom polynomial only generalizes naturally to {\em connected} cohomology theories ($h^{<0}=0$), as can be seen from the spectral sequence approach of Kazarian. Also a complex algebraic subvariety has a natural Poincar\'e dual in $K$-theory. However in other cohomology theories as cobordism theory you need extra information like a resolution. Some of these Poincar\'e duals when $G:=GL(n;\R)\times GL(n+k;\R)$ acts on $\hom(R^n,\R^{n+k})$ was calculated by Damon \cite{damon-residue} and Hayden \cite{hayden}. These calculations are quite difficult. But to calculate the whole avoiding ideal is usually much easier. Let us mention one---unpublished---result in this direction:
\begin{theorem}
  For any complex oriented cohomology theory $E$ Proposition \ref{pr-ideal} still holds if $c$ and $c'$ denote the corresponding Chern classes of $E$ and $\langle$\ $\rangle$ denotes the generated submodule over the coefficient ring (the $E$-cohomology of the point).
\end{theorem}
The proof is based on the simple fact that the equations defining the avoiding ideal are the same for any complex oriented cohomology theory.

\section{Classical results}\label{early}

\subsection{Singularities of smooth maps}
Consider a complex analytic map $f$ between complex analytic manifolds $N^n$ and $P^p$. Also let us fix a ``singularity'' $\eta$ (see a discussion below). An often occurring problem is the study of the set $\eta(f)$ of points $x\in N$ where $f$ has singularity $\eta$. Thom proved in \cite{thom} that the cohomology class represented by the closure of $\eta(f)$ in the cohomology ring of $N$ is equal to the value of a multivariable universal polynomial---depending only on $\eta$---when we substitute the characteristic classes of $N$ and the pull-backs of  the characteristic classes of $P$. In this way, if we know the polynomial, and the homotopy class of $f$, then we can tell the (co)homology class of $\overline{\eta(f)}$, obtaining direct geometric or topological consequences.

Unfortunately---for quite long---not many of these polynomials were known explicitly. Some known examples included works of Thom \cite{thom}, Porteous \cite{port}, Ronga \cite{rongaij} and Gaffney \cite{gaffney}. For a fuller list of references see \cite{avgl} or \cite{thompol}. In \cite{thompol} the second author applied (basically) the method of the present paper and found an algorithm to compute these polynomials.

In section \ref{singularities} we will show that these polynomials are Thom polynomials in the sense of the preceding section, and show that the method of
\cite{thompol} is exactly the application of Theorem \ref{alleqs}. In the meanwhile we will find what the good definition is for ``singularities''.

Let us remark that there is a parallel theory starting with real smooth manifolds and smooth maps between them. The method of the above mentioned authors gave results for the complex and the real case simultaneously. Considering more difficult singularities the two cases become essentially different, see also Subsection \ref{integer}.

While the techniques of the present paper seem to be powerful enough to reproduce most of the earlier results, let us mention that we never succeeded to re-calculate the results in \cite{rongaij}, i.e. the Thom polynomials associated to {\em all} second order Thom-Boardman singularities. The difficulty is that unlike the first order ones the higher order Thom-Boardman singularities are not orbits in the corresponding jet space and results of Section \ref{rest-equations}. cannot be directly applied.

\subsection{Lagrange and Legendre singularities}

A version of singularity theory of maps is obtained when we consider maps which come from some differential geometric situations. Usually the occurring maps are Lagrange maps between symplectic manifolds. Their singularities are called Lagrange singularities and their Thom polynomial theory has been studied and solved by Vassiliev and Kazarian, see \cite{vassiliev}, \cite{avgl}, \cite{kazarian95}, \cite{kazarian-ll}.

\subsection{Degeneracy loci.} Thom polynomials for group actions coming from algebraic geometry are usually called ``degeneracy loci formulas''. A review of the known formulas, as well as other enumerative properties of degeneracy loci, and discussions on Thom polynomials is Chow groups are found in \cite{fp}, see also \cite[Ch.~14]{fulton}. These investigations include the Giambelli-Porteous-Thom formula (see Section \ref{porteous} in the present paper), the interpretation of Schubert calculus (see also section \ref{ss}) and degeneracy loci formulas for other classical groups (see section \ref{lambda2}, and the original references \cite{jlp}, \cite{harris-tu} as well as generalizations by Kazarian \cite{kazarian-sym}). A recent result is an algorithm finding the degeneracy loci formulas for quiver representations associated with $A_n$ graphs, found by Fulton and Buch \cite{buch-fulton}, see more in section \ref{quiver}.

\subsection{Circle bundles} A spectacular topological example of Thom polynomials is given by Kazarian in \cite{kazarian-circle}. Here the Thom polynomial problem is the following. Let the (orientation preserving) diffeomorphism group of $S^1$ times that of $\R^1$ act on the space of maps from $S^1$ to $\R^1$ the natural way. Then in \cite{kazarian-circle} the Thom polynomials of a collection of orbits are computed. Translating this back to differential topology one obtains geometric interpretations of (multiples of) the powers of the Chern class of a complex line bundle.

\medskip

In the next sections we show how most of these results can be obtained using the method of Section \ref{rest-equations} together with some new results.

\section{Thom polynomials for $GL(n)$} \label{gl}

The naive approach would be calculating the Thom polynomials for irreducible representations and then finding the rules of calculating the Thom polynomials for the direct sums of representations if the Thom polynomials for the factors are known. However the orbit structure of a direct sum is usually quite complicated and typically has infinitely many orbits. To get a feel of the intricacies of the geometry of the direct sum let us have a look at the ``simplest'' example (for details see \cite{ss}):

\begin{example} $GL(n)$ acts on $V=\bigoplus_{i=1}^k\C^n\iso \hom(\C^k,\C^n)$. The kernel of the map $\phi\in V$ is the only invariant of the  $GL(n)$-action so the orbit space  $V/GL(n)$ will be the union of Grassmannians $\gr_i(\C^k)$. Using a cell decomposition for these ``moduli spaces'' we can define Thom polynomials which are going to be Schur polynomials.
\end{example}

 An extensive list of representations with finitely many orbits can be found in \cite{richardson} and \cite{kac}. A simple condition is $\dim G\geq \dim V$, otherwise we won't have an open orbit. The classification of irreducible representations of $GL(n)$ is well known together with a simple formula  for their dimensions (the ``hook length formula see e.g. \cite[p.50]{fulton-harris}). Recalling that taking the dual or tensoring with the one dimensional representations $\det^k$ doesn't change the dimension of a representation we get the following list:

\begin{proposition}
The irreducible representations $V$ of $GL(n)$ satisfying the condition $\dim GL(n)$ $\geq \dim V$ are the trivial representations $\C$, the standard representations $\C^n$, the symmetric and antisymmetric two forms  $S^2(\C^n)$, $\Lambda^2(\C^n)$, the adjoint representations of $SL(n)$, the representations $\Lambda^3(\C^6)$, $\Lambda^3(\C^7)$, $\Lambda^3(\C^8)$, and their duals and tensor products with $\det^k$.
\end{proposition}

The adjoint representations have infinitely many orbits (see a detailed discussion below) but all the others have a finite orbit structure. For $S^2(\C^n)$, $\Lambda^2(\C^n)$ this is the Sylvester theorem. For the representations $\Lambda^3(\C^6)$, $\Lambda^3(\C^7)$, $\Lambda^3(\C^8)$ see e.g. \cite[p.358]{fulton-harris}. We calculated some Thom polynomials for the exceptional cases in \cite{fnr}.

Taking the dual or tensoring with  $\det^k$ can change the Thom polynomials but in a controllable way so we concentrate on the list above.

The standard representation $\C^n$ has two orbits $\eta=\C^n\setminus\{0\}$ and $\{0\}$. Almost by definition the Thom polynomials are $\tp(\eta)=1$ and $\tp(0)=c_n$.

Orbits of $S^2(\C^n)$ and $\Lambda^2(\C^n)$ are determined by the rank and their Thom polynomials were calculated in \cite{jlp},\cite{harris-tu} (let $\Delta_{\lambda}$ denote the Schur polynomial associated with the partition $\lambda$, as in e.g. \cite{fp}):

\begin{theorem}[\cite{jlp},\cite{harris-tu}]\
\begin{itemize}
\item $\tp($orbit of $S^2(\C^n)$ with corank $r)=2^{r}\Delta_{(r,r-1,r-2,\ldots,1)}$.
\item $\tp($orbit of $\Lambda^2(\C^n)$ with corank $r)=\Delta_{(r-1,r-2,r-3,\ldots,1)}$.
\end{itemize}
\end{theorem}

Since these representations satisfy the Euler condition we can apply the method of restriction equations. On the case of $\Lambda^2(\C^n)$ let us demonstrate how simple these calculations are.

\subsection{Thom polynomials for $\Lambda^2(\C^n)$---restriction equations} \label{lambda2}

The second antisymmetric ($\Lambda^2$) power of the standard representation of $GL_n(\C)$ has finitely many orbits: the corank $r$ determines the orbit $\Sigma^r$---where $r=n,n-2,\ldots$. We can choose a representative of $\Sigma^r$ as
$$ M_r=\begin{pmatrix} 0 & I_{(n-r)/2} & 0 \\ -I_{(n-r)/2} & 0 & 0 \\ 0 & 0 & 0 \end{pmatrix}_{n\times n}.$$

Now we show how to obtain $\tp(\Sigma^r)=\Delta_{(r-1,r-2,\ldots,1)}(c)$ using our methods. (The case of $S^2(\C^n)$ is very similar.)

\begin{proof}
Easy computation shows that the maximal compact symmetry group of $M_r$ is $G_r=U(r) \times Sp(\frac{n-r}{2})$, with cohomology ring $H^*(BG_r)=\Z[c_1,\ldots, c_r, p_1,\ldots,p_{\frac{n-r}{2}}]$  ($c_i$ are the Chern classes of rank $2i$ and $p_i$ are the Pontryagin classes of rank $4i$). In other words this ring is $\Z[c_i,p_i]$, $i\in \Z$ with the ``substitutions''
$$c_{\hbox{negative}}=0,\qquad p_{\hbox{negative}}=0,\qquad c_{>r}=0,\qquad p_{>\frac{n-r}{2}}=0.$$
Let us call these substitutions $(*r)$. With this notation the map between cohomology groups induced by the inclusion $j_r:G_r\to GL_n(\C)$ is
$$j_r^*:c_i \mapsto \sum_{l=0}^{\infty} c_{i-2l}p_l\Big|_{(*r)},$$
and the Euler class of $\Sigma^r$ is
$$\Delta_{(r-1,r-2,\ldots,1,0,0,0,\ldots)}(c)\Big|_{(*r)}.$$
To prove the theorem let us compute $j_s^*(\Delta_{(r-1,r-2,\ldots,1)}(c))$. It is
$$\Delta_{(r-1,r-2,\ldots,1,0,0,\ldots)}\Big(\sum_{l=0}^{\infty} c_{i-2l}p_l|_{(*s)}\Big)=$$
$$=\sum_{f:\N\to\N} p_1^{|f^{-1}(1)|}p_2^{|f^{-1}(2)|}\ldots \cdot \Delta_{(r-1-2f(1),r-2-2f(2),\ldots)}(c)\Big|_{(*s)}.$$
The $\Delta$'s in the last sum corresponding to $f\not\equiv 0$ are all zero, since the $i$\textsuperscript{th} and the $i+f(i)$\textsuperscript{th} rows are identical for the greatest $i$ such that $f(i)\not=0$. The term corresponding to $f\equiv 0$ is $e_r$ if $r=s$, and has only 0's in the first row if $s<r$. These prove the principal and the homogeneous equations, so the theorem.
\end{proof}

We have calculated the whole avoiding ideals for these orbits in \cite{fnr}.

\subsection{The adjoint representations}
We show that the adjoint representation  of $GL(n)$ has no interesting Thom polynomials. This is true for the adjoint representation of a large class of Lie groups including the complex semisimple  Lie groups but the proof is essentially the same so we only indicate how the general case works.

The adjoint representation $V=\ad(GL(n))$ is not irreducible: $V=\ad(SL(n))\oplus \C$. We state the following proposition for  $V$ but it is true for $\ad(SL(n))$, too and the proof works the same way.

\begin{proposition} If the invariant subset $\xi\subset V$ has a Thom polynomial then $\tp(\xi)=0$ or 1. \end{proposition}

\begin{proof}
The orbits of $V$ are described by the Jordan normal forms. A generic $v\in V$ has $n$ different eigenvalues and the stabilizer group $G_\eta$ of the orbit $\eta=GL(n)v$ is isomorphic to $GL(1)^n$. The inclusion $G_\eta\to GL(n)$ induces a map $j^*_\eta: H^*_{GL(n)}(pt)\to H^*_{GL(1)^n}(pt)$. This map is injective. (This fact is usually called the splitting lemma. For semisimple  Lie groups the corresponding map is still injective rationally by a theorem of Borel.) In other words $\ob_\eta=0$. It implies that for any invariant subset $\xi\subset V$ for which $\eta\not\subset \bar\xi$ the avoiding ideal $\av_\xi \subset\ob_\eta=0$ so  $\tp(\xi)=0$. And if $\bar\xi$ contains every generic orbit then $\bar\xi=V$ and $\tp(\xi)=1$.
\end{proof}

\begin{remark}
In other words the first obstruction for finding a section  of an $\ad(GL(n))$-bundle with different eigenvalues is zero. But there are other obstructions: Suppose that $E$ is a complex vector bundle over a simply connected base and $\ad(E)$ is its adjoint bundle. If $\ad(E)$ admits a  section with different eigenvalues at every point then $E=\bigoplus L_i$ where the $L_i$'s are the one dimensional eigenspace bundles (if the base is not simply connected we get only an $n$-line distribution), but not all bundles split.
\end{remark}

\subsection{Representations of $GL(2)$}
The relevant irreducible representations of $GL(2)$ are of the form
$$S^n\C^2=\{\mbox{homogeneous polynomials of order $n$ in two variables}\}.$$
As we mentioned at the beginning of the section Thom polynomials for the other irreducible representations can be easily obtained from these. For $n\geq 4$ these representations  have families of orbits. A homogeneous polynomial $p=p(x,y)$ defines a 0-dimensional variety $V(p)$ in the projective line $\P^1$. Let us call the points of $V(p)$---counted with multiplicity---the {\em roots} of $p$. Since the cross ratio of 4 points  is an invariant of the $GL(2)$-action we will get families of orbits for $n\geq 4$. However if $p$ has at most 3 different roots then the orbit of $p$  doesn't have a family. It turns out that the existence of these orbits allows us to calculate the Thom polynomial for the orbits  with at most two different roots:
$$\eta_i:=\mbox{orbit of }\ x^iy^{n-i} \ \ \mbox{for} \ \ i=0,\dotsc,[n/2].$$
Notice that $x^iy^{n-i}$ is in the same orbit as $x^{n-i}y^i$. This is the reason why we get a slightly different formula for $n$ even and $n$ odd.

First we calculate the geometric input: The stabilizers (more precisely the maximal compact subgroups of the stabilizers) $G_{\eta_i}$  of these orbits and the the action of $G_{\eta_0}$ on the normal space $N_{\eta_0}$ of $\eta_0$:

\begin{proposition} \mbox{} \label{geomdata}
\begin{enumerate}[(i)]
\item $G_{\eta_i}$ is the image of the homomorphism $h_i:U(1)\to GL(2)$ where $h_i(\alpha)=\smx{\alpha^{n-i}}{}{}{\alpha^i}$.
\item  $N_{\eta_0}=\langle x^n,x^{n-1}y,\dotsc,x^2y^{n-2}\rangle$.
\end{enumerate}
\end{proposition}
The calculation are elementary so we omit the proof.

Proposition \ref{geomdata} implies that
$$h_i^*(c_1)=(n-2i)c_1 \ \ h_i^*(c_2)=-i(n-i)c_1^2 \ \ \mbox{for} \ \ i=1,\dotsc,[n/2] \ \ \mbox{and} $$
$$h_0^*(c_1)=c_1  \ \ h_0^*(c_2)=0  \ \ \mbox{and} \ \ e(\eta_0)=n!c_1^{n-1}.$$

Using Theorem \ref{alleqs} we see that if $\eta_i$ is not in the closure of an orbit $\eta$ then $\tp(\eta)$ is divisible by $i(n-i)c_1^2+(n-2i)^2c_2$. Taking the principal equation into account we get that for example:

\begin{theorem}$$\twocase{\tp(\eta_0)=}{\phantom{\frac n2c_1}n\prod\limits_{i=1}^t \big(i(n-i)c_1^2+(n-2i)^2 c_2 \big)}{n \ \ \hbox{\rm is odd}}{n\frac n2c_1\prod\limits_{i=1}^t\big(i(n-i)c_1^2+(n-2i)^2 c_2 \big)}{n \ \ \hbox{\rm is even}.}$$
\end{theorem}

We can associate an invariant subset $(n_1,\dotsc,n_k)$ to any partition of $n$ where the partition encodes the required multiplicities of the roots of the polynomials in the subset. We calculated the Thom polynomials  $\tp(n_1,\dotsc,n_k)$ with different methods in \cite{fnr-coincident}. Recently B. K\H om\H uves found a closed formula using incidences in the sense of \cite{thompol}.

\begin{remark} The Thom polynomials $\tp(1^k,n-k)$ were calculated in  \cite{kirwan} by Kirwan for the $SL(2)$-action (notice that for the $SL(2)$-action $c_1=0$ so the $SL(2)$-Thom polynomial contains less information).
\end{remark}

\section{The classical case: Giambelli-Thom-Porteous formula} \label{porteous}

In this section we show how to recover the classical Thom polynomial formula (the so called Giambelli-Thom-Porteous formula, see \cite{thom}, \cite{port}) in our theory. We choose our field to be $\C$ and our cohomology theory to be $H^*(\ ,\Z)$. Suppose that $f:N\to P$ is a smooth map of manifolds. The  Giambelli-Thom-Porteous formula describes the  cohomology class defined by $\Sigma_s(f)$, the subset of $N$ where $df$ has corank~$s$.

In terms of the theory described above we calculate the Thom polynomials  of the representation $\rho=\hom\bigl(\rho(n),\rho(n+k)\bigr)$ of the group $G:=GL(n)\times GL(n+k)$ on the linear space $\C^{(n+k)\times n}$  where $\rho(n)$ is the standard representation of $GL(n)$. So $(R,L)\in G$ acts on an $(n+k)\times n$ matrix $X$ by: $(R,L)\cdot X:=LXR^{-1}$. We can assume that $k\geq0$.

As it is well known the orbits $\Sigma_s$ of this action are characterized by corank. A representative from $\Sigma_s$ is $X_s:=\begin{pmatrix} 0 & 0 \\ 0 & I_{n-s} \end{pmatrix}_{(n+k)\times n}$.  The maximal compact stabilizer subgroup of $X_s$ is
$$G_s:=G_{X_s}=\Big\{\ \Big( \bigl(\begin{smallmatrix} A & 0 \\ 0 & C \end{smallmatrix}\bigr), \bigl(\begin{smallmatrix} B & 0 \\ 0 & C \end{smallmatrix} \bigr)\Big) \ \Big|\ (A , B, C) \in U(s)\times U(s+k)\times U(n-s) \Big\}.$$

An invariant normal slice $\Sigma_s$ at $X_s$ is $N_s=\bigl\{\smx{M^{(s+k)\times s}}{0}{0}{0}\bigr\}$. It implies that the---complex---codimension of $\Sigma_s$ is $s(s+k)$. To determine the principal equation for the Thom polynomial associated to $\Sigma_s$ we need two data: $e(\Sigma_s)$, the $G$-equivariant Euler class of $\Sigma_s$, and the map $H^*(Bi_s): H^*(BG) \to H^*(BG_s)$ where $i_s$ is the inclusion of $G_s$ into $G$. We will use the notation

\begin{align*}
H^*(BG)&= \Z[R_1,\dotsc, R_n, L_1,\dotsc,L_{n+k}]\\
&=\Z^{S_n\times S_{n+k}} [r_1,\dotsc,r_n, l_1,\dotsc, l_{n+k}]\\
H^*(BG_s)&= \Z[A_1,\dotsc, A_s, B_1, \dotsc, B_{s+k}, C_1, \dotsc, C_{n-s}]\\
     &= \Z^{S_s\times S_{s+k}\times S_{n-s}}[a_{n-s+1},\dotsc, a_n, b_{n-s+1}, \dotsc, b_{n+k}, c_1, \dotsc, c_{n-s}],
\end{align*}

\noindent where the capitals mean universal Chern classes while the lower letters mean Chern roots, and  e.g. $\Z^{S_n\times S_{n+k}}[\ ]$ means the part of the polynomial ring $\Z[\ ]$ invariant under the action of $S_n\times S_{n+k}$, i. e. the permutations of the $a_i$'s and $b_j$'s.

Notice that the  $G$-equivariant Euler class of $\Sigma_s$ is equal to $e(N_s)$, the $G_s$-equivariant Euler class of the normal slice $N_s$. The action of $(A,B,C)\in G_s$ on   $N_s$ is given by changing $M$ to $BMA^{-1}$. So written in terms of Chern roots:
\[e(\Sigma_s)=e(N_s)=\prod_{i=n-s+1}^n \quad\prod_{j=n-s+1}^{n+k} (b_j-a_i)\in H^{s(s+k)}(BG_s).\]
Since $e(\Sigma_s)$ is not zero for any $s$ the representation satisfies the Euler condition and by Theorem \ref{unicity} the restriction equations have a unique solution. Even more is true:

\begin{proposition}\label{princ} For every orbit $\Sigma_s$ the principal equation $H^*(Bi_s)\tp(\Sigma_s)=e(N_s)$ has a unique solution. \end{proposition}
\begin{proof}
Since $G_s\subset G_{s+1}$ the restriction map $H^*(Bi_s)$ factors through $H^*(Bi_{s+1})$ so the homogeneous equations contain no extra information.
\end{proof}

\begin{remark}
  The same condition---that the principal equation has a unique solution---applies to the representations $\Lambda^2\C^n$ and $S^2\C^n$ studied in Section \ref{gl}. However in that cases it was easy to prove that the Thom polynomial satisfies the homogeneous equations.
\end{remark}
\begin{remark}
The fact $G_s\subset G_{s+1}$ has even stronger consequences. It implies that $\ob_{\Sigma_{s+1}}\subset\ob_{\Sigma_s}$ therefore the avoiding ideal
\begin{equation}\av_{\Sigma_{s+1}}=\bigcap_{i\leq s}\ob_{\Sigma_i}=\ob_{\Sigma_s}=\ker H^*(Bi_s).\label{avo}\end{equation}
 Using the explicit description of $H^*(Bi_s)$ below it makes it easy to check whether a characteristic class in $H^*(BG)$ belongs to $\av_{\Sigma_{s+1}}$. Also, since $\av_{\Sigma_{s+1}}\cap H^d(BG)=0$ for $d<\codim \Sigma_{s+1}$ formula (\ref{avo}) implies the injectivity of $H^{\codim \Sigma_s}(Bi_s)$ which in turn implies Proposition \ref{princ} directly without using the unicity theorem \ref{unicity}.
\end{remark}

The map $H^*(Bi_s)$  is given by (again in terms of Chern roots):

$$\begin{array}{lll}
   r_i \mapsto &&\left\{\begin{aligned} c_i &\quad&\text{if} &\quad& i\leq n-s \\ a_i &\quad&\text{if}
&\quad& i>n-s \end{aligned}\right. \end{array}\ \ \ \ \ \ \ \ \ \begin{array}{lll}
  l_i \mapsto &&\left\{\begin{aligned} c_i &\quad&\text{if} &\quad& i\leq n-s \\ b_i &\quad&\text{if} &\quad& i>n-s
\end{aligned}\right.  \end{array}$$

\begin{lemma} Using the notation
$$\frac{1+L_1t+L_2t^2+\ldots+L_{n+k}t^{n+k}}{1+R_1t+R_2t^2+\ldots+R_nt^n}=  1+H_1t+H_2t^2+\ldots$$
 $\det (H_{s+i-j})_{(s+k)\times (s+k)}$ maps to $\prod_{i=n-s+1}^n \prod_{j=n-s+1}^{n+k} (b_j-a_i)$ under the map $H^{s(s+k)}(Bi_s)$.
\end{lemma}

\begin{proof}
The image $H'_i$ of $H_i$ are  the coefficients of the Taylor series
$$\frac{\prod (1+b_it) \prod (1+c_kt)}{\prod (1+a_j t)\prod(1+c_kt)}=\frac{\prod (1+b_it)}{\prod (1+a_j t)}=
\frac{1+B_1t+B_2t^2+\ldots+B_{s+k}t^{s+k}}{1+A_1t+A_2t^2+\ldots+A_st^s}.$$

But $\det (H'_{s+i-j})_{(s+k)\times (s+k)}$ equals to the resultant of the two polynomials $1+B_1t+B_2t^2+\ldots+B_{s+k}t^{s+k}$ and $1+A_1t+A_2t^2+\ldots+A_st^s$. (This is a less known form of the resultant $R(p,q)$ which can be obtained by multiplying the Sylvester-matrix by a matrix obtained from the coefficients of the Taylor series of $1/p$ see \cite[p.87]{arba}.) On the other hand the resultant is equal to the product of differences of the roots of the two polynomials.
\end{proof}

This lemma  and Proposition \ref{princ} together proves that there is only one polynomial that satisfies the principal equation:
 $$\tp(\Sigma_s)=\det (H_{s+i-j})_{(s+k)\times (s+k)}.$$
Notice that a $\rho$-bundle in this case is a pair of vector bundles $E,F$ of rank $n$ and $n+k$. In the classical situation of a map $f:N\to P$ these bundles are $TN$ and $f^*TP$, and $H_i$ can be interpreted as the $i$\textsuperscript{th} Chern class of the virtual bundle $F\ominus E$. Also notice that the formula doesn't depend on $n$.

Using Section \ref{ideal} one could easily calculate the whole avoiding ideal $\av_{\Sigma_s}$, but here we only give the result, since it has been computed by P. Pragacz.
\begin{proposition}[\cite{pragacz-enum}, {\cite[Section 4.2]{fp}}] \label{pr-ideal} $\av_{\Sigma_s}=\{\Delta_\lambda(c,c'): \lambda\supset (s+k)^s\}.$
\end{proposition}

\section{Singularities}\label{singularities}

In this section we show how our theory applies to the case of singularities of maps between manifolds---the case where Thom polynomials were originally defined by Thom in \cite{thom}. We will work over the complex field, so manifolds and maps are assumed to be complex analytic. What we really show is that the equations we get by Theorem \ref{alleqs} for the Thom polynomials of simple singularities are the same as were studied and solved in \cite{thompol}, so we will not repeat their solution here.

Now we recall some standard definitions of singularity theory (see e.g. \cite{avgl}): $\Eo(n,n+k)$ will  be the vector  space of smooth germs $(\C^n,0)\to (\C^{n+k},0)$. We will think of $\Eo(n,n+k)$ as a subset of $\Eo(n+a,n+k+a)$ by trivial unfolding. Fixing $k$ let $\Eo(\infty,\infty+k)$ be the union (or formally the direct limit):  $\cup_{n=0}^{\infty} \Eo(n,n+k)$. This space  will play  the role  of $V$  of  the general theory. Also we have maps $u_\infty:\Eo(n,n+k)\to \Eo(\infty,\infty+k)$ ($u_\infty$ stands for infinite {\em unfolding}). Let $\hol(\C^n,0)$ denote the group of biholomorphism germs of $(\C^n,0)$. The group
$$\A(n,n+k):=\hol(\C^n,0)\times\hol(\C^{n+k},0)$$
acts  on $\Eo(n,n+k)$ by $(\phi,\psi)\cdot f:=\psi\circ f \circ \phi^{-1}$. Similarly the limit group
$$\A(\infty,\infty+k):=\cup_{n=0}^{\infty}\hol(\C^n,0)\times\hol(\C^{n+k},0)$$
acts on $\Eo(\infty,\infty+k)$ by the same formula.

We will mainly be concerned with the bigger---{\em contact}---groups
$$\K(n,n+k):=\{(\varphi,M):\varphi \in \hol(\C^n,0), M\ \text{is a germ}\ (\C^n,0)\to\hol(\C^{n+k},0) \},$$
acting on $\Eo(n,n+k)$ by $\rho^{\K(n,n+k)}(f)= M(x)\circ f \circ \phi^{-1}(x)$. The limit group $\K(\infty,\infty+k)$ acts on $\Eo(\infty,\infty+k)$ by the same formula. This group will play the role of $G$ in the general theory. We use the notations $\Eo,\A$ and $\K$ if the value of $n$ ($n=\infty$ allowed) is clear from the context.

So, consider the action of $\K(\infty,\infty+k)$ on $\Eo(\infty,\infty+k)$. The nicest orbits are the so called {\em simple} ones: an orbit (or a representative) is simple, if a neighborhood intersects only finitely many different orbits. Simple orbits will be strata in an appropriate Vassiliev stratification. Let $\eta$ be a simple orbit, and let us choose a representative $f\in \Eo(n,n+k)$ with  minimal $n$. In other words we choose  a  minimal dimensional  representative  with  $\eta=$ orbit  of $u_{\infty}(f)$. Such an $f$ (defined up to $\K(n,n+k)$-equivalence)  is called a {\em  genotype} for $\eta$ in  \cite[p.~157]{avgl}. The contact automorphism group  $\stab^\K(f)=\{(\phi,M)\in \K(n,n+k)| (\phi,M)\cdot f =f\}$ and the analogously defined $\stab^\A(f)$ are not finite dimensional (moreover they do not possess convenient topologies) so we need the following definition---inspired by the classical Bochner theorem---using that $GL(n)\times GL(n+k)\subset\A(n,n+k) \subset\K(n,n+k)$.

\begin{definition}[{\cite{janich}}] \label{mcdef} If $M$ is a subgroup of $\A(n,n+k)$ or $\K(n,n+k)$ then $M$ is {\em compact} if  $M$ is conjugate to a compact subgroup $N \subset GL(n)\times GL(n+k)$.
\end{definition}

Luckily enough the groups $\stab^\K(f)$ and $\stab^\A(f)$ share many properties with finite dimensional groups, as follows.

\begin{theorem}\label{maxcompact}
  \mbox{}
\begin{enumerate}
 \item $\stab^\K(f)$ ( $\stab^\A(f)$) has a maximal compact subgroup $G_f=G^\K_f$ ($G^\A_f$). \label{emax}
 \item Any two maximal compact subgroups are conjugate. \label{aconj}
 \item $B\stab^\K(f)\htpy BG^\K_f$ and $B\stab^\A(f)\htpy BG^\A_f$. \label{contr}
\end{enumerate}
\end{theorem}
The proof of (\ref{emax}) and  (\ref{aconj}) can be found in \cite{janich}, \cite{wall2}. As we mentioned $\stab(f)$ does not possess convenient topology, so strictly speaking  $B\stab(f)$ is not defined. However, it is possible to define the notion of $\stab(f)$-principal bundle over a {\em smooth manifold} and $BG_f$ classifies those bundles (\cite[Thm 1.3.6]{rrphd} or  \cite{rl}). So from our point of view we can replace $B\stab(f)$ by $BG_f$. In particular we have $B\K(n,n+k)\htpy B\A(n,n+k)\htpy BGL(n)\times BGL(n+k)$.

\begin{remark}
Theorem \ref{maxcompact} allows us to use the same algorithm to calculate the Thom polynomials as in the finite dimensional case.

Definition \ref{mcdef} and Theorem \ref{maxcompact} shows that by choosing $f$ carefully from its $\K$-orbit, we can assume that $G^\K_f\subset GL(n)\times GL(n+k) $, so we have representations $\mu_0$, $\mu_1$  of $G^\K_f$ on the source space $V_0$ and target space $V_1$ respectively. By part  \ref{aconj} of Theorem \ref{maxcompact} the isomorphism classes of these representations are uniquely  defined. The groups  $G^\K_f$ and representations $\mu_0$, $\mu_1$ were calculated for low codimensional singularities in \cite{rl}.
\end{remark}

For the identification of $G_\eta$ for $G=\K(\infty,\infty+k)$ we cannot directly use Theorem \ref{maxcompact}, but it is not difficult to get around:

\begin{definition}
$$G_\eta:=G_f\times U(\infty)\subset \stab(\eta)$$
\noindent
where the inclusion of $ U(\infty)$ into $\stab(\eta)$ corresponds to the diagonal action on the unfolding dimensions.
\end{definition}

Though $G_\eta$ is not compact in any reasonable sense, it is still true that $B\stab(\eta)\htpy BG_\eta$ (in the sense of our remark after Theorem \ref{maxcompact}) and as we will see the $ U(\infty)$ summand acts trivially on $N_\eta$ anyway.

Below we explain how to calculate the two inputs of the algorithm for computing the Thom polynomials for an orbit $\eta$, i. e. the map $ H^*(Bi): H^*(BG)\to  H^*(BG_\eta)$ and the representation $\rho_\eta: G_\eta\to GL(N_\eta)$.

\begin{proposition} \label{cohom}
The homomorphism $H^*(B\K)=H^*(BU(\infty)\times U(\infty))\to H^*(BG_{\eta})$ induced by the inclusion $G_{\eta}\to \A\subset \K$ is given by
$$\begin{array}{ccc}
\Z[\u{a},\u{b}] & \a & H^*(BG_f) \otimes \Z[\u{d}] \\
a & \mapsto & c(\mu_0)\cdot d \\
b & \mapsto & c(\mu_1)\cdot d, \\
\end{array}$$
where $\u{a}=a_1,a_2,\ldots$ and $\u{b}=b_1,b_2,\ldots$ are the universal Chern classes of the two factors of $U(\infty)\times U(\infty)$; $a=1+a_1+a_2+\ldots$, $b=1+b_1+b_2+\ldots$ are the total Chern classes, and the definitions for $\u{d}$, $d$ are similar. The class $c(\mu_i)$ is the total Chern class of the vector bundle $E_{\mu_i}=EG_f\times_{\mu_i}V_i$ over $BG_f$. \qed
\end{proposition}

Let us turn to our second goal. It is enough to calculate the $G_f^\K$-action $\rho_f$ on the normal space $N_f^\K$ to the $\K(n,n+k)$-orbit in $\Eo(n,n+k)$ since $u_{\infty}(N_f^{\K})$ is normal to $\eta$ in $\Eo(\infty,\infty+k)$ too, and  once the action is $G_f$-invariant, it is also $G_{\eta}$-invariant with the  trivial action of the $U(\infty)$ factor.

The representation  $\rho_f$ is also explicitly computable for low codimensional singularities: The miniversal unfolding $F$ of $f$ is a stable germ
$(\C^n\oplus U,0) \to (\C^{n+k}\oplus U,0)$ where $U$ is the unfolding space. The  group $G_f^\K$ acts linearly on $U$: let us denote this  representation by $\mu_U$. Then we have $\rho_f\iso\mu_0\oplus\mu_U$. We can also see that the source space of $F$ is naturally isomorphic to $N_f^\K$. (For more details see \cite{wall2} or \cite{rl}.)

\begin{remark}
Calculation of the \tpol s for the stable orbits of the $\A(\infty,\infty+k)$-action doesn't give anything new: every such orbit is a dense open subset of a $\K(\infty,\infty+k)$-orbit---we usually use the same notation for the two orbits---so their \tpol s are the same.
\end{remark}

Now, we might as well write down the equations of Theorem \ref{alleqs} for the \tpol s for $\K$, but it is possible to simplify these equations as we will see in Theorem \ref{eqs}.

We can also calculate \tpol s for $\K(n,n+k)$. These results are not independent: If $\eta$ is an orbit of $\K(n,n+k)$ then its $d$-dimensional unfolding  $u_d(\eta)$ is an orbit of $\K(n+d,n+k+d)$ with the same codimension. This is a consequence of the fact that the codimension of the $\K$-orbit can be read off from its local algebra, which doesn't change by trivial unfolding. To understand the connection between the \tpol\ of $\eta$ and of $u_d(\eta)$ we look at the unfolding map $u_d: \K(n,n+k) \to \K(n+d,n+k+d)$. It induces a map
 \begin{multline*}
H^*(Bu_d):H^*(B  \K(n+d,n+k+d)\iso \Z[a_1,\dotsc,a_{n+d},b_1,\dotsc,b_{n+k+d}] \to \\ H^*(B \K(n,n+k))\iso \Z[a_1,\dotsc,a_{n},b_1,\dotsc,b_{n+k}]
   \end{multline*}

\noindent ($d=\infty$ is allowed) such that

$$ H^*(Bu_d)(a_i)=\left\{\begin{aligned} a_i &\quad&\text{if} &\quad& i\leq n \\ 0 &\quad&\text{if} &\quad& i>n \end{aligned}\right. \qquad\qquad H^*(Bu_d)(b_i)=\left\{\begin{aligned} b_i &\quad&\text{if} &\quad& i\leq n+k \\ 0 &\quad&\text{if} &\quad& i>n+k \end{aligned}\right. $$

\begin{proposition}\label{tp-unf}$ \tp(\eta)=H^*(Bu_d)\big(\tp(u_d(\eta))\big)$ where $u_d$ denotes the d-fold trivial unfolding ($d=\infty$ is allowed).
 \end{proposition}

This is a special case of a more general fact used frequently in calculating Thom polynomials:

\begin{proposition}  \label{pullback}
Let $V_1$ be a $G_1$-vector space and $V_2$ be a $G_2$-vector space Let $\phi:G_1\to G_2$ be a homomorphism and $J:V_1\to V_2$ a $\phi$-equivariant map transversal to the $G_2$-action. Suppose that $\eta\subset V_2$ has a Thom polynomial $\tp(\eta)$. Then $\tp(\phi^{-1}(\eta))=\phi^*\tp(\eta)$.
\end{proposition}

\begin{proof} This is a straightforward generalization of the pullback property of the ordinary Poincar\'e dual.  \end{proof}

\begin{proof}[Proof of Proposition \ref{tp-unf}] The unfolding map $u_d$ is transversal so we can apply Proposition \ref{pullback}. \end{proof}

So using the generators $a_i, b_i$ for all of these groups we can say that for $d$ large enough ($n+d\geq \codim\eta$) unfoldings don't change the \tpol. Pursuing these arguments further, we find another important property of these \tpol s (first appeared probably in \cite{damonphd}):

\begin{proposition}[{folklore, see \cite{avgl}}] \label{quotient} $\tp(\eta) \in Q$ where $Q$ is the  subring of $H^*(B\K)$ generated by $1,h_1,h_2,\ldots$, where $1+h_1+h_2+\ldots=\frac{1+b_1+b_2+\ldots}{1+a_1+a_2+\ldots}$.
\end{proposition}

Before getting into the proof we need some definitions:
\begin{definition}
If $A^n$ and $B^{n+k}$ are vector bundles over a manifold $M$, then let $\Eo(A,B)$ denote the $\Eo (n,n+k)$-bundle over $M$ such that
$$\Eo_m(A,B)=\{\text{\rm germs of smooth maps }(A_m,0)\to (B_m,0)\}.$$
\end{definition}

Using Theorem \ref{maxcompact} it is  not difficult to see that
$$a_i(\Eo(A,B))=c_i(A)\ \ \text{and }\ b_i(\Eo(A,B))=c_i(B)$$
through the obvious identifications. (In fact it also follows from  Theorem \ref{maxcompact} that every $\Eo (n,n+k)$-bundle is isomorphic to a bundle of the form $\Eo(A,B)$ but we don't use this in this paper.)

\begin{definition}
Let $C^d$ be a $d$-dimensional vector bundle. Then
$$u_C: \Eo(A,B)\to \Eo(A\oplus C,B\oplus C)$$
denotes the {\em twisted unfolding map}:
$$u_C(\varphi):=\varphi\oplus\id_C.$$
\end{definition}

{\it Proof of Proposition \ref{quotient}.}
We have a \cd
$$\xymatrix{ E_{\rho^{\K(n,n+k)}}\ar[r]^{\tilde u_d}  & E_{\rho^{\K(n+d,n+k+d)}} \\
\Eo(A,B) \ar[r]^-{u_c} \ar[u]^{\tilde k}   & \Eo(A\oplus C,B\oplus C)  \ar[u]^{\tilde k_C}
}$$
\noindent where $\tilde k$ and $\tilde k_C$ are the bundle maps induced by the corresponding classifying maps. It shows that
$$k^*\tp(\eta)=k^*_C\tp(u_d(\eta))$$
Choose $C$ to be the `inverse' of $A$. Then using that we can think of $\tp$ as a polynomial of $\u a$ and $\u b$ we get (with some abuse of notation):
$$\tp(\eta)\big(c(A),c(B)\big)=\tp(u_d(\eta))\big(1,c(B)/c(A)\big)$$
\qed

\begin{definition}
Let $\tpq(\eta)$ be the unique polynomial with the property
$$\tpq(\eta)(1,h_1,\dotsc)=\tp(\eta)(\u a,\u b),$$
\noindent  where
$1+h_1+h_2+\ldots=\frac{1+b_1+b_2+\ldots}{1+a_1+a_2+\ldots}$.
\end{definition}

Proposition \ref{cohom} is enough to write down the equations for the Thom polynomials of simple singularities in the cohomology ring of $B\A$ or $B\K$, namely in $\Z[\u{a}, \u{b}]$. In the light of Proposition \ref{quotient} we write it in terms of the `quotient' variables $h_i$:

\begin{theorem} \label{eqs}
$$\tpq(\eta)(c(\theta)) =\begin{cases}
\hbox{\rm Euler class of } E_{\rho_\eta} & \hbox{\rm if  } \theta=\eta  \qquad \qquad \hbox{`principal equation'}\cr 0
& \hbox{\rm if  } \theta\not> \eta \qquad\qquad  \hbox{`homogeneous equations'}, \end{cases}$$
\noindent where $c(\theta)=c(E_{\mu_1(\theta)})/c(E_{\mu_0(\theta)})$.
\end{theorem}

These are exactly the equations that were solved for many cases in \cite{thompol}.

\begin{remark} \label{rel-euler}
This representation of $\K$ doesn't satisfy the Euler condition but a closer look at the proof of the unicity theorem \ref{unicity} shows that the restriction equations have a unique solution in the codimension range where there are only finitely many orbits and the Euler classes of the normal spaces are non zero.
\end{remark}

\begin{remark}
At this point we would like  to comment on the history of these ideas. The systematic study  of classifying spaces of the  symmetry groups of singularities and a  powerful construction  out of these spaces was pioneered by A. Sz\H ucs (see e.g. \cite{szucs1}). In the language of the present paper he calculated $G_\eta$ for several singularities, and described a way how to glue these spaces together to get a space whose algebraic topological properties can be translated to differential topological theorems. He applied his construction to various differential topological questions, such as e.g. the cobordism groups of maps  with given singularities  (e.g. \cite{szucs2}, \cite{szucs4}, \cite{szucs5}, \cite{szucs-morin}).  A general method  of calculating more of the symmetry groups was given in \cite{rrphd}, \cite{rrsz}, \cite{rl}. In the present paper we explored the fact that roughly speaking (the `source space' in) Sz\H ucs' construction is a union  of strata in the stratification of $B\A$ defined  by the $\A$-action, and that the ideas fruitful there (e.g. Thom polynomial calculations  \cite{thompol}) turn out to be fruitful viewing any $G$-action.
\end{remark}

\begin{example} \label{1,1}
Let us start with $V=\Eo(1,k+1)$. The group $\K(1,k+1)$ acts on it as usual. Here the finite codimensional strata will be the contact orbits $A_n^1(k)$ represented by the germ $x\mapsto (x^{n+1},0,\dotsc,0)$. Just like above one can write up the equations for the Thom polynomials of these strata. Carrying out the computation (or using Proposition \ref{pullback}) one finds that here the `principal' and `homogeneous' equations are  {\em enough} (see Remark \ref{rel-euler}) to determine the Thom polynomials:
\begin{equation}\tp(A_n^1(k))=\prod_{j=0}^k\prod_{i=1}^n (b_j-ia)\in H^*\big(BU(1)\times BU(k+1)\big)=\Z^{S_{k+1}}[a,b_0,\dotsc,b_k].\label{tp-ai}\end{equation}

Since the unfolding map $u_\infty:\Eo(1,k+1)\to\Eo(\infty,k+\infty)$ is transversal to the orbits by Proposition \ref{pullback} we have $u_\infty^*\big(\tp(A_n(k))\big)=\tp(A_n^1(k))$. The homomorphism $u_\infty^d$ is injective for $d\leq 2k+2$ so this way we get a simple way to calculate the Thom polynomials of $A_2(k)$ (the so called Ronga formulas). $u_\infty^{3k+3}$ has a kernel so formula (\ref{tp-ai}) is not enough to calculate the Thom polynomials of $A_3(k)$. (It makes it easier however. We published a closed formula for $\tp(A_3(k))$ in \cite{bfr}.)
\end{example}

\section{Other applications of the method of RE}\label{applications}

\subsection{Thom polynomials for singularities with integer coefficients} \label{integer}
There is a theory of {\em real} singularities parallel to the complex case discussed in Section \ref{singularities}. References for the rich geometry and topology of this local theory are e.g. \cite{avgl} and \cite{wdp}. The globalization of the theory, i.e. in the study of the Thom polynomials of real singularities has two levels: one can ask for the Thom polynomials with mod 2 coefficients or with integer coefficients. The mod 2 case can be basically solved by the method of \cite{b-h}. Already in this case, and very essentially in the $\Z$-coefficient case one has to answer the question: what (formal linear connections of singularities) defines a Thom polynomial? The answer uses a detailed analysis of the Kazarian spectral sequence (or at least its 0\textsuperscript{th} row, the so called {\em Vassiliev universal complex}). Both the answer to this question and the actual Thom polynomials are presented in \cite{int}.

If a singularity $\eta$ is not cooriented then there is no Thom polynomial with integer coefficients. However the avoiding ideal may not be empty. A simple way to find element in $\mathcal A_\eta$ is to look at the ``realization" of an element in $\mathcal A_{\eta_\C}$. For example $p_1$ is an obstruction for finding an immersion of an $n$-dimensional manifold into $\R^{n-1}$ since $\tp(\Sigma^1_\C(-1))=c_2$. $p_1$ as an obstruction was used in \cite{sst}.

\subsection{Quivers}\label{quiver}
A surprisingly rich class of representations are the quiver representations. Let $E_0$ be the set of vertices and $E_1$ be the set of edges of a connected oriented graph (double edges, loops allowed). The tail and the head of an edge $e$ is denoted be $t(e)$ and $h(e)$. If there is a vector space $V_v$ assigned to any vertex $v$ then we can consider the group $G=\prod_{v\in E_0} GL(V_v)$ and its action on the vector space
$$\bigoplus_{e\in E_1} \hbox{Hom}(V_{t(e)},V_{h(e)})$$
given by
$$\Big(M_v\Big)_{v\in E_0}\ \cdot \Big(\phi_e\Big)_{e\in E_1}= \Big(M_{h(e)}\circ \phi_e\circ M^{-1}_{t(e)}\Big)_{e\in E_1}.$$
This representation is called the quiver representation associated to the graph and the dimension vector $($dim$(V_v))_{v\in E_0}$.

This representation, including its orbit structure was thoroughly studied in representation theory (see e.g. \cite{ars}). It turns out that this action has finitely many orbits if and only if the graph is of Dynkin type (with some orientations on the edges).

The Thom polynomials of the action associated to the $A_n$ series were studied by Buch and Fulton \cite{buch-fulton}. They showed on one hand that these Thom polynomials generalize many objects in algebraic combinatorics (see also \cite{fp}), including different versions of Schubert polynomials (for which there has been no explicit determinantal formula known, they are usually computed by recursion). On the other hand in \cite{buch-fulton} an algorithm is given to compute the Thom polynomials of any orbit of $A_n$. Buch and Fulton conjecture a formula for these Thom polynomials, and prove it for special cases.

In \cite{fr-quiver} the authors applied the method of the present paper to quiver representations associated with arbitrary Dynkin graphs, and thus obtained a straightforward procedure (but not a formula) to get any Thom polynomials.

\subsection{Schur and  Schubert polynomials}\label{ss}
 The cohomology ring structure of Grassmann and flag manifolds are governed by Schur and Schubert polynomials, see e.g. \cite{fp}. In a recent application (\cite{ss}) the authors found a way to obtain these Schur and Schubert polynomials as Thom polynomials.

The starting point is that we act on the vector space $\Hom(\C^n,\C^p)$ by triangular matrices from both sides (flag case) or by triangular matrices from one side and $GL$ from the other side (Grassmann case). Then the orbits will correspond to Schubert varieties, whose equivariant Poincar\'e dual are the  double Schubert (flag case) and double Schur (Grassmann case) polynomials. If we disregard one set of indeterminates we obtain the ordinary Schubert and ordinary Schur polynomials.

So it is enough to apply the method of the present paper to the above triangular and half-triangular actions, and we recover (or give a new definition as well as a new way to compute) double and simple Schur and Schubert polynomials. We proved a Jacobi-Trudi type determinantal formula for the double Schur polynomials and gave an illuminating proof for the Lascoux-Sch\"utzenberger recursion (\cite{ls}) for double Schubert polynomials.

\section{Projective Thom polynomials and degree calculations}\label{sec:ptp}

If a group $G$ acts linearly on a vector space $V$ then this action $\rho$ induces an action $\p\rho$ of $G$ on the projective space $\p V$. If the image of $\rho$ contains the scalars then there is a bijection between the orbits of $\p\rho$ and the non zero orbits of $\rho$. Strangely enough the projective Thom polynomials formally contain more informations than the affine ones: Suppose that $\eta$ is a $d$ complex codimensional invariant subset of $\rho$ and $\p\eta$ is the corresponding invariant subset of $\p\rho$. Then $\tp(\p\eta)$ is an element in $H^{2d}_G(\p V)\iso H^{2d}(BG)[\xi]/\prod \xi-\beta_i$ where $\beta_i$ are the weights of the representation $\rho$. So $\tp(\p\eta)=\sum p_i\xi^i$ where $p_i\in H^{2(d-i)}(BG)$. It is easy to see that $p_0$ is the affine Thom polynomial $\tp(\eta)$ and $p_d$ is the degree of the closure of $\p\eta$:
\begin{equation}
  \label{eq:coeffs}
  p_0=\tp(\eta) \qquad \qquad p_d=\deg(\eta)
\end{equation}
In fact the main application of the projective Thom polynomial is that we can calculate the degree of certain varieties.

The main result of this section is that a simple substitution into the affine Thom polynomial $\tp(\eta)$ provides the projective Thom polynomial. To state the result we need to give names to the generators of $H^*(BG)$. Let $m:U(1)^n\to G$ be a (coordinatized) maximal torus of $G$ and let $\alpha_i$ are the corresponding roots. So by the Borel theorem (or splitting principle) $\tp(\eta)$ is a polynomial in the roots $\alpha_i$ and $\tp(\p\eta)$ is a polynomial in the roots $\alpha_i$ and $\xi$. We assume that the image of $\rho$ contains the scalars i.e. there is a homomorphism $\phi:GL(1)\to G$ and a non zero integer $q$ such that $\rho\circ\phi(\lambda)=\lambda^qv$ for all $v\in V$, $ \lambda\in GL(1)$. We assume that $\im m\supset \im \phi|_{U(1)}$ so we have a homomorphism $\tilde\phi:U(1)\to U(1)^n$ such that $\phi|_{U(1)}=m\circ\tilde\phi$. The  homomorphism $\tilde\phi$ is necessarily of the form $\tilde\phi(t)=(t^{w_1},\dotsc,t^{w_n})$ where $t\in U(1)$ and $w_i$ are integers. Notice that the choice of $\phi$ is not unique.

\begin{theorem}[{\cite{fnr}}]\label{ptp} Let $\rho:G\to GL(V)$ be a representation of the Lie group $G$ such that the image of $\rho$ contains the scalars. Let $\alpha_i$, $q$, $w_i$ be as above and let $\eta$ be an invariant subset of $\rho$. Then
$$\tp(\p\eta)(\alpha_1,\dotsc,\alpha_n,\xi)=\tp(\eta)(\alpha_1+\frac{w_1}q\xi,\dotsc,\alpha_n+\frac{w_n}q\xi).$$
\end{theorem}

\begin{corollary} \label{degree} Using the notation of Theorem \ref{ptp}
$$\deg(\p\eta)=q^{-d}\tp(\eta)(w_1,\dotsc,w_n).$$
\end{corollary}
This is a generalization of results on degrees of certain degeneracy loci of Porteous \cite{port}, Harris-Tu \cite{harris-tu}, Fulton \cite{fulton}. In \cite{fnr} this formula is used to calculate  the degree of the dual of some of the Grassmannians $Gr_k(\C^n)$.

\section{The \kss}\label{kss}

The Kazarian  spectral sequence gives interesting relations (e.g. linear equations or bounds) for the number of different strata---i.e. in  some  sense, the  number of  different singularity types---in a fixed codimension. In this section we are giving two examples of these.

Let us briefly recall the Kazarian spectral sequence \cite{kazarian}. Associated with the group representation $\rho:G\to GL(V)$ we consider the universal $\rho$-bundle: $BV:=EG\times_{\rho}V\longrightarrow BG$ (the letter B in $BV$ stands for {\em Borel} construction). A key observation is that if $\eta$ is an invariant subset of $\rho$ then $\eta$ can be identified in every fiber. Let their collection be $B\eta$ ($=$ the Borel construction with $\eta$). So a stratification of $V$ by invariant submanifolds induces a stratification of $BV$. If $F_i$ is the union of $B\eta$'s with codim $\eta$ at most $i$, then the spectral sequence associated with the filtration $\emptyset\subset F_0\subset F_1\subset\ldots\subset BV$ is called the Kazarian spectral sequence. Now suppose that the stratification we started with is special: it satisfies Vassiliev's condition \cite[8.6.5]{vassiliev}---i.e. it is locally finite, the stabilizer subgroups in one stratum are constant and the moduli spaces stratum$/G$ are contractible. For example if the representation has finitely many orbits and the stratification is the orbit stratification then Vassiliev's condition is trivially satisfied. Then using excision and Thom isomorphism one can easily see that the $E_1$ term of this spectral sequence is the following: the $q$\textsuperscript{th} column contains the---possibly twisted---cohomology groups of the classifying spaces of the stabilizer subgroups of $q$-codimensional strata. In this general situation the Thom polynomial is defined as the edge homomorphism of this spectral sequence: $E^{p,0}_2\to E^{p,0}_{\infty}\subset H^*(BV)=H^*(BG)$.

So this spectral sequence clearly organizes the points to consider when there are neighboring dimensional strata and one has to glue some together to obtain one which satisfies the conditions we put in Section \ref{tp}. Considering complex representations this spectral sequence often (but not always) degenerates at $E_1$. Since this was the situation in almost all the examples we considered, we restricted ourselves to the simple version of Section~\ref{obstruction}.

The study of this spectral sequence was initiated in \cite{ab} to study the betti numbers of certain moduli spaces.
Let us now, however, give an  other and simpler application of the Kazarian spectral sequence.

Consider the \kss\ associated to the representation of $GL(n)\times GL(n+k)$ (over $\C$) discussed in Section \ref{porteous}. It collapses at $E_1$ since the odd rows and columns are zero.  Particularly interesting  is the limiting case $n=\infty$. We assume for simplicity that $k=0$. One finds that in the
$s$\textsuperscript{th} column one has to write the  cohomology groups of $U(s)^2\times U(\infty)$, see Section \ref{porteous}. Summing up the ranks in the skew diagonals we must obtain the ranks of the cohomologies of $U(\infty)^2$. Some combinatorics shows that ``one can drop a $\times U(\infty)$ term'' from everywhere, i.e. one can write the cohomology groups of $U(s)^2$ in the $s$\textsuperscript{th} column and get the ranks of the cohomologies of $U(\infty)$ (i.e. $1,1,2,3,5,7,11,15,\dotsc$) in the skew diagonals, as follows:

$$\begin{array}{c|cccccccccc}
&   &    &   &   &    &   &  &  &  & \\
8. & . &  5 & . & . & 14 & . & . & .  & . & 16\\
6. & . &  4 & . & . &  8 & . & . & .  & . & 10\\
4. & . &  3 & . & . &  5 & . & . & .  & . & 5\\
2. & . &  2 & . & . &  2 & . & . & .  & . & 2\\
0. & 1 &  1 & . & . &  1 & . & . & .  & . & 1\\
 \cline{2-11}
 \multicolumn{2}{r}{ 0.}& 2. & 4.& 6.& 8. & 10. & 12.  & 14. & 16. & 18. \\
\end{array}$$

In this table only the ranks of the free Abelian groups are written (i.e. $a$ is written instead of $\Z^a$) and only the terms with two even coordinates are indicated, since everything else is 0. This leads to the combinatorial identity
    $$\pi(n,\![1,1,\dotsc])\! \!=\!\pi(n-1,\![2,0,0,\dotsc]) +\pi(n-4,\![2,2,0,0,\dotsc])+\pi(n-9,\![2,2,2,0,0,\dotsc])+\dotsb\!,$$
where $\pi(n,[a_1,a_2,\dotsc])$ denotes the number of degree $n$ monomials in terms of $a_i$ copies of variables of degree $i$. This identity---already known to Euler as a very effective way to compute the number of partitions---can be directly proved by using Young diagrams (thanks to A. Blokhuis for these informations).

Let us now turn to the Kazarian spectral sequence of $\Eo(\infty,k+\infty)$ of Section \ref{singularities}. For simplicity let $k=0$. The list of simple singularities (for codimension $\leq 8$) is given in the following table.
$$\begin{array}{r|ccccccccc}
\hbox{codim}_{\C} & 0 & 1 & 2 & 3 & 4 & 5 & 6 & 7 & 8 \\
\hline
  & A_0 & A_1 & A_2 & A_3 & A_4 & A_5 & A_6 & A_7 &  A_8 \\
& & & & & I_{2,2} & I_{2,3} & I_{2,4} & I_{2,5}   &  I_{2,6} \\
& & & & &         &         & I_{3,3} & I_{3,4}   &  I_{3,5} \\
& & & & &         &         &         &           &  I_{4,4} \\
& & & & &         &         &         & (x^2,y^3) &  (x^2+y^3,xy^2) \\
\end{array}$$
The maximal compact symmetry groups of these singularities
can be computed as in \cite{rl}, which gives us the $E_1$ term for the Kazarian spectral sequence. Let us suppose however, that we only know the easy part of the classification, i.e. only up to codimension 3. The symmetry group of the $A_i$ singularity is easily computed as $U(1)\times U(\infty)$ for $i>0$. So hereby we show the (degenerated) spectral sequence with the `$U(\infty)$ terms dropped' (as above), and writing only the ranks of the occurring groups:
$$\begin{array}{r|rrrrr}
10.&0 &1 &1 &1 &  \\
8.&0 &1 &1 &1 & \\
6.&0 &1 &1 &1 &  \\
4.&0 &1 &1 &1 &  \\
2.&0 &1 &1 &1 &  \\
0.&1 &1 &1 &1 &n \\
\cline{2-6}
\multicolumn{2}{r}{ 0.} & 2. &4. &6. &8.  \\
\end{array}$$

Here thus $n$ is the number of singularities of complex codimension 4 (which we assumed not to know). The point here is that the value of $n$ can be found from the spectral sequence above, by observing that the the sums of skew diagonals should be equal to the ranks of the cohomologies of $BU(\infty)$: $1,1,2,3,5,\ldots$ (the number of partitions). This gives us that $n$ must be 2. So we could predict the number of different strata of codimension $d$ knowing only information about strata of codimension $<d$. The interested reader can extend the above spectral sequence using the further classification of singularities, and try to speculate about the number of different strata in codimensions $16,17,18,\ldots$ (but be careful about moduli spaces!).

Applications of this method in algebraic geometry and singularity theory are subject to further study.

\bibliography{sing}
\bibliographystyle{alpha}

\end{document}